\documentclass[11pt,reqno]{amsart}
\usepackage{amssymb,amscd,amsbsy,mathrsfs}
\renewcommand{\a}{\alpha}

\renewcommand{\d}{\delta}

\newcommand{\z}{\zeta}

\renewcommand{\l}{\lambda}

\newcommand{\s}{\sigma}
\renewcommand{\t}{\tau}

\newcommand{\f}{\varphi}
\renewcommand{\o}{\omega}

\renewcommand{\L}{\Lambda}

\renewcommand{\O}{\Omega}

\newcommand{\E}{{\mathscr E}}

\newcommand{\F}{{\mathscr F}}

\newcommand{\I}{{\mathscr I}}

\newcommand{\1}{{\bf 1}}

\newcommand{\C}{{\Bbb C}}
\newcommand{\T}{{\Bbb T}}

\newcommand{\dd}{{\Bbb D}}
\newcommand{\R}{{\Bbb R}}
\newcommand{\Z}{{\Bbb Z}}

\newcommand{\0}{{\boldsymbol{0}}}

\newcommand{\bs}{\boldsymbol}

\newcommand{\bS}{{\boldsymbol S}}

\newcommand{\rf}[1]{(\ref{#1})}

\newcommand{\df}{\stackrel{\mathrm{def}}{=}}
\newcommand{\dist}{\operatorname{dist}}
\newcommand{\Ker}{\operatorname{Ker}}

\newcommand{\supp}{\operatorname{supp}}

\newcommand{\const}{\operatorname{const}}

\newcommand{\eeq}{\end{equation}}
\newcommand{\beq}{\begin{equation}}
\newcommand{\bay}{\begin{eqnarray}}
\newcommand{\ba}{\begin{align*}}
\newcommand{\ea}{\end{align*}}
\newcommand{\ey}{\end{eqnarray}}
\newcommand{\bey}{\begin{eqnarray*}}
\newcommand{\eey}{\end{eqnarray*}}

\newcommand{\be}{\infty}

\newcommand{\bl}{\blacksquare}

\newcommand{\Pf}{{\bf Proof. }}

\def\mathR{\mbox{$\bf \rm I\makebox [1.2ex] [r] {R}$}}

\newtheorem{thm}{\hspace{\parindent}Theorem}[section]

\newtheorem{cor}[thm]{\hspace{\parindent}Corollary}
\newtheorem{lem}[thm]{\hspace{\parindent}Lemma}

\usepackage{amsmath}

\begin{document}

\theoremstyle{remark}
\newtheorem{rem}[thm]{Remark}
\newtheorem*{rem*}{Remark}

\numberwithin{equation}{section}

\numberwithin{equation}{section}

\title{Functions of noncommuting operators under perturbation of class $\bs{\bS_p}$}
\author{A.B. Aleksandrov and V.V. Peller}
\thanks{The research of the first author is partially supported by RFBR grant 17-01-00607.
The publication was prepared with the support of the
RUDN University Program 5-100}
\thanks{Corresponding author: V.V. Peller; email: peller@math.msu.edu}


\newcommand\ite{{\Li\,\hat\otimes_{\rm i}\,L^\be}}
\newcommand\te{{\Li\,\hat\otimes\,L^\be}}
\newcommand\cE{\mathcal{E}}
\newcommand\up{\upsilon}
\newcommand\Li{{\rm Lip}}
\newcommand\fM{\frak M}
\newcommand\cZ{\mathcal{Z}}
\newcommand\dg{\frak D}
\newcommand\pt{\hat\otimes_{\rm i}}
\newcommand{\fI}{{\frak I}}
\newcommand{\mt}{{\mathcal T}}
\newcommand\mC{\mathcal{C}}
\newcommand\mQ{\mathcal{Q}}
\newcommand\mR{\mathcal{R}}
\newcommand\Wp{\stackrel{1}{W}}
\newcommand\Wv{\stackrel{2}{W}}

\newcommand\sG{\mathscr G}

\newcommand{\qm}{\quad\mbox{and}\quad}

\newcommand\mX{\mathcal{X}}
\newcommand\mY{\mathcal{Y}}
\newcommand\mB{\mathcal{B}}
\newcommand\fn{\frak n}
\newcommand\bn{\bs{\fn}}
\newcommand\fm{\frak m}
\newcommand\fp{\frak p}
\newcommand\mZ{\mathcal Z}

\newcommand\mI{\mathcal{I}}

\newcommand{\OL}{{\rm OL}}
\newcommand{\COL}{{\rm CL}}
\newcommand{\com}{{\,\rm com}}
\newcommand{\SA}{{\bf SA}}
\newcommand{\No}{{\bf N}}
\newcommand{\Un}{{\bf U}}
\newcommand{\Pro}{{\bf P}}
\newcommand{\USA}{{\bf USA}}
\newcommand{\fF}{{\frak F}}
\newcommand{\CO}{{\bf C}}
\newcommand{\sgn}{\operatorname{sgn}}
\newcommand{\card}{\operatorname{card}}
\newcommand{\Lip}{\operatorname{Lip}}
\newcommand{\Ba}{\big(B_{\be,1}^1\big)_+(\T^2)}
\newcommand\CA{{\rm C}_{\rm A}}

\newcommand{\BL}{{\rm BOL}}

\newcommand\J{\mathscr{J}}

\def\mathR{\mbox{$\bf \rm I\makebox [1.2ex] [r] {R}$}}

\begin{abstract}
In this article we prove that for $p>2$, there exist pairs of self-adjoint operators $(A_1,B_1)$ and $(A_2,B_2)$ and a function $f$ on the real line in the homogeneous Besov class $B_{\be,1}^1(\R^2)$ such that the differences
$A_2-A_1$ and $B_2-B_1$ belong to the Schatten--von Neumann class $\bS_p$ but
$f(A_2,B_2)-f(A_1,B_1)\not\in\bS_p$. A similar result holds for functions of contractions. We also obtain an analog of this result in the case of triples of self-adjoint operators for any $p\ge1$.
\end{abstract}

\maketitle

\setcounter{section}{0}
\section{\bf Introduction}
\setcounter{equation}{0}
\label{In}

\

The paper is devoted to the study of Schatten--von Neumann properties of the increments $f(A_2,B_2)-f(A_1,B_1)$ under perturbation of pairs of not necessarily commuting bounded self-adjoint operators.

It was Farforovskaya who discovered in \cite{F} that for Lipschitz functions $f$ on the real line $\R$ the condition that $A-B\in\bS_1$ for self-adjoint operators $A$ and $B$ does not imply that $f(A)-f(B)\in\bS_1$. Here $\bS_1$ stands for trace class, see \cite{GK}.

Later it was shown in \cite{Pe1} and \cite{Pe2} for functions $f$ in the homogeneous Besov space $B_{\be,1}^1(\R)$ (see \S\:\ref{prel}) the condition that $A-B\in\bS_1$ does imply that 
$f(A)-f(B)\in\bS_1$ and
\bay
\label{OLS1}
\|f(A)-f(B)\|_{\bS_1}\le\const\|f\|_{B_{\be,1}^1}\|A-B\|_{\bS_1}.
\ey
Moreover, it was shown there that under the same assumption this inequality also holds for all Schatten--von Neumann classes $\bS_p$, $1\le p\le\be$:
\bay
\label{OLSp}
\|f(A)-f(B)\|_{\bS_p}\le\const\|f\|_{B_{\be,1}^1}\|A-B\|_{\bS_p}.
\ey
Note that for $p=\be$, by $\|\cdot\|_{\bS_p}$ we mean   operator norm. 

It is well known (see \cite{AP1}) that $f$ satisfies \rf{OLS1} for arbitrary self-adjoint operators with trace class difference if and only if $f$ is an opetrator Lipschitz function, i.e.,
$$
\|f(A)-f(B)\|\le\const\|f\|_{B_{\be,1}^1}\|A-B\|
$$
for arbitrary self-adjoint operators $A$ and $B$. We refer the reader to \cite{AP1} for detailed information on operator Lipschitz functions.

It was shown in \cite{PS} that for $p\in(1,\be)$, there exists a positive number $c_p$ such that
\bay
\label{PS}
\|f(A)-f(B)\|_{\bS_p}\le c_p\|f\|_{\rm Lip}\|A-B\|_{\bS_p}
\ey
whenever $A$ and $B$ are self-adjoint operators with $A-B\in\bS_p$.

Analogs of \rf{OLSp} for functions of $n$-tuples of self-adjoint operators were obtained in \cite{APPS} for $n=2$ and in \cite{NP} for 
$n>2$. We would also like to mention the paper \cite{KPSS}, in which an analog of 
\rf{PS} for $n$-tuples of self-adjoint operators was obtained.

Functions of pairs of noncommuting self-adjoint operators can be defined with the help of double operator integrals, see \cite{ANP}. We discuss this issue in detail in \S\:\ref{Fnko}. In particular, the functions $f(A,B)$ can be defined for arbitrary bounded self-adjoint operators $A$ and $B$ and for every $f$ in the Besov class $B_{\be,1}^1(\R^2)$, see \S\:\ref{Fnko}.

It was shown in \cite{ANP} that for $p\in[1,2]$, the following estimate of Lipschitz type holds in the norm of the Schatten--von Neumann ideal  $\bS_p$:
\bay
\label{Litip}
\|f(A_2,B_2)-f(A_1,B_1)\|_{\bS_p}\le
\const\|f\|_{B_{\be,1}^1(\R^2)}\max\{\|A_2-A_1\|_{\bS_p},\|B_2-B_1\|_{\bS_p}\}
\ey
for arbitrary pairs $(A_1,B_1)$ and $(A_2,B_2)$ of not necessarily commuting self-adjoint operators such that $A_2-A_1\in\bS_p$ and $B_2-B_1\in\bS_p$.
Moreover, an analog of this result for functions of unitary operators was also obtained there.

However, in the same paper \cite{ANP} it was established that for $p>2$, there are no such Lipschitz type estimates as well as there are no such estimates in the operator norm. 

Nevertheless, the authors of \cite{ANP} failed to answer the important question of whether for $p\in(2,\be)$, there exist a function $f$ of class $B_{\be,1}^1(\R^2)$ and pairs $(A_1,B_1)$ and $(A_2,B_2)$ of self-adjoint operators such that $A_2-A_1\in\bS_p$ and $B_2-B_1\in\bS_p$ but $f(A_2,B_2)-f(A_1,B_1)\not\in\bS_p$.

The main result of this paper (Theoremа \ref{ne vkhodit00sa}) gives an affirmative answer to this question.

In Section \ref{szhali!} of this article we consider a similar problem for functions of not necessarily commuting contractions. Recall that in \cite{AP} the following analog of inequality \rf{Litip} was obtained for $p\in[1,2]$:
\bay
\label{pin[1,2]}
\|f(T_2,R_2)-f(T_1,R_1)\|_{\bS_p}\le
\const\|f\|_{(B_{\be,1}^1)_+(\T^2)}\max\{\|T_2-T_1\|_{\bS_p},\|R_2-R_1\|_{\bS_p}\}
\ey
for an arbitrary function $f$ in the Besov class $\Ba$ (see \S\;\ref{prel}) of analytic functions and for arbitrary pairs of contractions $(T_1,R_1)$ and $(T_2,R_2)$ on Hilbert space satisfying $T_2-T_1\in\bS_p$ and $R_2-R_1\in\bS_p$. Recall that an operator $T$ is called a {\it contraction} if $\|T\|\le1$. 

Note that it was shown in \cite{Pe5} that inequality \rf{pin[1,2]} holds for $p\in[1,\be]$ in the case of pairs of {\it commuting contractions} $(T_1,R_1)$ and $(T_2,R_2)$ on Hilbert space.

It was shown in \cite{AP} that such an inequality of Lipschitz type does not hold for $p>2$ as well as in the operator norm.

The main result of Section \ref{szhali!} of this paper is that for $p\in(2,\be)$, there exist a function $f$ of class $\Ba$ and pairs $(T_1,R_1)$ and $(T_2,R_2)$ of contractions on Hilbert space such that
$T_2-T_1\in\bS_p$ and $R_2-R_1\in\bS_p$ but $f(T_2,R_2)-f(T_1,R_1)\not\in\bS_p$.
Moreover, one can even construct such pairs of unitary operators $(T_1,R_1)$ and $(T_2,R_2)$.

Finally, in \S\:\ref{triple} we obtain analogs of the results of \S\:\ref{samosop2} for functions of triples of arbitrary bounded self-adjoint operators.

\

\section{\bf Besov spaces}
\setcounter{equation}{0}
\label{prel}

\

In this article we deal with the homogeneous Besov class $B_{\be,1}^s(\R^2)$, $s=1,2$,
of functions on the real line and with the Besov class $\Ba$ of functions on the two-dimensional torus that are analytic on the bidisk. We refer the reader to the book \cite{Pee} and the papers \cite{ANP} and \cite{AP} for information on these classes.
Let $w$ be an infinitely differentiable function on $\R$ such
that
\bay
\label{w}
w\ge0,\quad\supp w\subset\left[\frac12,2\right],\quad\mbox{and} \quad w(t)=1-w\left(\frac t2\right)\quad\mbox{for}\quad t\in[1,2].
\ey

Consider the functions $W_n$, $n\in\Z$, on $\R^d$ such that 
$$
\big(\F W_n\big)(x)=w\left(\frac{\|x\|_2}{2^n}\right),\quad n\in\Z, \quad x=(x_1,\cdots,x_d),
\quad\|x\|_2\df\left(\sum_{j=1}^dx_j^2\right)^{1/2}.
$$

With each tempered distribution $f\in{\mathscr S}^\prime\big(\R^d\big)$, we
associate the sequence $\{f_n\}_{n\in\Z}$,
\bay
\label{fn}
f_n\df f*W_n.
\ey
The formal series
$
\sum_{n\in\Z}f_n
$
is a Littlewood--Paley type expansion of $f$. 
This series does not necessarily converge to $f$. 

Initially we define the (homogeneous) Besov class $\dot B^1_{\be,1}\big(\R^d\big)$ as the space of 
$f\in{\mathscr S}^\prime(\R^n)$
such that
\bay
\label{<be}
\|f\|_{B^1_{\be,1}}\df\sum_{n\in\Z}2^n\|f_n\|_{L^\be}<\be.
\ey
According to this definition, the space $\dot B^1_{\be,1}(\R^n)$ contains all polynomials
and all polynomials $f$ satisfy the equality $\|f\|_{B^1_{\be,1}}=0$. Moreover, the distribution $f$ is determined by the sequence $\{f_n\}_{n\in\Z}$
uniquely up to a polynomial. It is easy to see that the series 
$\sum_{n\ge0}f_n$ converges in ${\mathscr S}^\prime(\R^d)$.
However, the series $\sum_{n<0}f_n$ can diverge in general. Obviously, the series
\bay
\label{ryad}
\sum_{n<0}\frac{\partial f_n}{\partial x_j},\quad
1\le j\le d,
\ey
converges uniformly on $\R^d$. 
Now we say that  $f$ belongs to the homogeneous Besov class $B^1_{\be,1}(\R^d)$ if \rf{<be} holds and 
$$
\frac{\partial f}{\partial x_j}=\sum_{n\in\Z}\frac{\partial f_n}{\partial x_j},\quad
1\le j\le d.
$$

A function $f$ is determined uniquely by the sequence $\{f_n\}_{n\in\Z}$ up
to a a constant, and a polynomial $g$ belongs to 
$B^1_{\be,1}\big(\R^d\big)$
if and only if it is a constant.

In the definition of $B_{\be,1}^2(\R^d)$ we should replace inequality \rf{<be} with the following one:
$$
\|f\|_{B^2_{\be,1}}\df\sum_{n\in\Z}2^{2n}\|f_n\|_{L^\be}<\be
$$
and condition \rf{ryad} with a similar condition that involves not only first order partial derivatives, but also second order.

\medskip

Studying periodic functions on $\R^d$ is equivalent to studying functions on the $d$-dimensional torus $\T^d$. To define Besov spaces on $\T^d$, we consider a function $w$ satisfying \rf{w} and define the trigonometric polynomials $W_n$, $n\ge0$, by
$$
W_n(\z)\df\sum_{j\in\Z^d}w\left(\frac{|j|}{2^n}\right)\z^j,\quad n\ge1,
\quad W_0(\z)\df\sum_{\{j:|j|\le1\}}\z^j,
$$
where 
$$
\z=(\z_1,\cdots,\z_d)\in\T^d,\quad j=(j_1,\cdots,j_d),\quad\mbox{and}\quad
|j|=\big(|j_1|^2+\cdots+|j_d|^2\big)^{1/2}.
$$
For a distribution $f$ on $\T^d$ we put
$$
f_n=f*W_n,\quad n\ge0,
$$
and we say that $f$ belongs the Besov class $B_{\be,1}^1(\T^d)$ if
$$
\sum_{n\ge0}2^n\|f_n\|_{L^\be}<\be.
$$
Note that locally the Besov space $B_{p,q}^s(\R^d)$ coincides with the Besov space
$B_{p,q}^s$ of periodic functions on $\R^d$.

We also define the Besov class $\big(B_{\be,1}^1\big)_+(\T^d)$ as the subspace of
$B_{\be,1}^1(\T^d)$ of holomorphic functions on the polydisk $\dd^d$, i.e.,
$$
\big(B_{\be,1}^1\big)_+(\T^d)=\big\{f\in B_{\be,1}^1(\T^d):\widehat f(j_1,j_2,\cdots,j_d)=0~{\rm whenever}~\min\{j_1,j_2,\cdots,j_d\}<0\big\}.
$$

We refer the reader to \cite{Pee} and \cite{AP1} for more detailed information on Besov spaces.

\

\section{\bf Functions of noncommuting operators operators}
\label{Fnko}

\

In this section we describe briefly how to define functions of nocommuting operators.

Let $A$ and $B$ be not necessarily commuting bounded self-adjoint operators on Hilbert space. Functions $f(A,B)$ are defined by the following double operator integral
$$
f(A,B)=\iint_{\R\times\R}f(x,y)\,dE_A(x)\,dE_B(y),
$$
where $E_A$ and $E_B$ are the spectral measures of $A$ and $B$, whenever the double operator integral on the right makes sense. Note that if the spectra $\s(A)$ and 
$\s(B)$ are contained in closed intervals $\I$ and $\J$, then $f(A,B)$ depends only on the restriction of $f$ to $\I\times\J$.

Recall that double operator integrals, i.e., expressions of the form 
\bay
\label{DOI}
\iint_{\R\times\R}\Phi(x,y)\,dE_1(x)Q\,dE_1(y),
\ey
appeared first in \cite{DK}. Here $\Phi$ is a measurable function, $E_1$ and $E_2$ are spectral measures and $Q$ is a bounded linear operator. Later Birman and Solomyak developed in \cite{BS1}--\cite{BS3} a beautiful theory of double operator integrals. For the double operator integral \rf{DOI} to be defined, the function $f$ has to satisfy certain assumptions.

The maximal class of functions $\Phi$, for which the double operator integral \rf{DOI} can be defined is called the class of {\it Schur multipliers} with respect to the spectral measures $E_1$ and $E_2$. There are several characterizations of the class of Schur multipliers, see \cite{Pe1}. 

We are going to use the following: $\Phi$ is a Schur multiplier with respect $E_1$ and $E_2$ if and only if $\Phi$ belongs to the {\it integral projective tensor product} $L^\be_{E_1}\widehat\otimes_{\rm i}L^\be_{E_2}$ of the $L^\be$ spaces 
$L^\be_{E_1}$ and $L^\be_{E_2}$, i.e., $\Phi$ admits a representation
$$
\Phi(x,y)=\int_\O\f(\o,x)\psi(\o,y)\,d\s(\o),
$$
where $(\O,\s)$ is a $\s$-finite measure space, and $\f$ and $\psi$ are measurable
functions such that
$$
\int_\O\|\f(\o,\cdot)\|_{L^\be_{E_1}}\|\psi(\o,\cdot)\|_{L^\be_{E_2}}\,d\s(\o)<\be.
$$
In this case 
$$
\iint_{\R\times\R}\Phi(x,y)\,dE_1(x)Q\,dE_1(y)=
\int_\O\left(\int\f(\o,x)\,dE_1(x)\right)Q\left(\int\psi(\o,y)\,dE_2(y)\right).
$$
In particular, $\Phi$ is a Schur multiplier if $\Phi$ belongs to the {\it projective tensor product} $L^\be_{E_1}\widehat\otimes L^\be_{E_2}$, i.e., $\Phi$ admits a representation
$$
\Phi(x,y)=\sum_{n\ge1}\f_n(x)\psi_n(y),
$$
where $\f_n$ and $\psi_n$ are measurable functions such that
$$
\sum_{n\ge1}\|\f_n\|_{L^\be_{E_1}}\|\psi_n\|_{L^\be_{E_2}}<\be.
$$

Let us also mention that $\Phi$ is a Schur multiplier if and only if it belongs to the Haagerup tensor product $L^\be_{E_1}\otimes_{\rm h}L^\be_{E_2}$, see
\cite{AP1}, \cite{Pe+} and \cite{Pe1}.

Let us return to functions of pairs $A$ and $B$ of noncommuting self-adjoint operators. Suppose that $f\in L^\be(\R^2)$ and its Fourier transform $\F f$ belongs to $L^1(\R^2)$.
Then 
$$
f(x,y)=\iint_{\R\times\R}(\F f)(s,t)e^{{\rm i}sx}e^{{\rm i}ty}\,ds\,dt,
$$
and so $f$ is a Schur multiplier with respect to arbitrary Borel measures on $\R$ and 
the operator $f(A,B)$ can be defined by
$$
f(A,B)=\iint_{\R\times\R}(\F f)(s,t)
e^{{\rm i}sA}e^{{\rm i}tB}\,ds\,dt.
$$

In particular, we can define $f(A,B)$ for functions $f$ in the homogeneous Besov class $B_{\be,1}^1(\R^2)$. Indeed, it suffices to show that such a function locally coincides with a function whose Fourier transform is in $L^1(\R)$. This in turn can be reduced to the following fact for periodic functions: if 
$f\in B_{\be,1}^1(\T^2)$, then 
\bay
\label{absskho}
\sum_{j,k\in\Z}\big|\hat f(j,k)\big|<\be.
\ey
Indeed, 
\begin{align}
\label{absolyutno!}
\sum_{j,k\in\Z}\big|\hat f(j,k)\big|
\le\sum_{n\ge0}\sum_{j,k\in\Z}\big|\hat f_n(j,k)\big|
&\le\const\sum_{n\ge0}2^n
\left(\sum_{j,k=-2^{n+1}}^{2^{n+1}}\big|\hat f_n(j,k)\big|^2\right)^{1/2}
\nonumber\\[.2cm]
&\le\const\sum_{n\ge0}2^n\|f_n\|_{L^\be}\le\const\|f\|_{B_{\be,1}^1}.
\end{align}

There is another way to define $f(A,B)$ for functions $f$ in $B_{\be,1}^1(\R^2)$, see \cite{ANP}. Without loss of generality we may assume that $\|A\|<\pi$, 
$\|B\|<\pi$ and $f$ is a $2\pi$-periodic function of Besov class $B_{\be,1}^1$.
Then we can represent $f$ as an element of the projective tensor product
$C(\T)\hat\otimes C(\T)$ in the following way:
$$
f(s,t)=\sum_{n\ge1}\sum_{j=-2^n}^{2^n}e^{{\rm i} sj}
\left(\sum_{j=-2^n}^{2^n}\widehat f_{n-1}(j,k)e^{{\rm i} tk}
\right),\quad s,~t\in\R.
$$
Clearly, this representation allows us to estimate the tensor norm of $f$:
$$
\|f\|_{C(\T)\hat\otimes C(\T)}\le
\sum_{n\ge1}\big(2^{n+1}+1\big)\|f_{n-1}\|_{L^\be(\T^2)}\le\const\|f\|_{B_{\be,1}^1}
$$
by \rf{<be}. Under these assumptions, the function $f(A,B)$ can be defined by
$$
f(A,B)=\sum_{n\ge1}\sum_{j=-2^n}^{2^n}e^{{\rm i}jA}
\left(\sum_{j=-2^n}^{2^n}\widehat f_{n-1}(j,k)e^{{\rm i}kB}
\right).
$$

Note here that in the case when $A$ and $B$ have finite spectra, we can define
$f(A,B)$ for arbitrary functions $f$ on $\s(A)\times\s(B)$:
$$
f(A,B)=\sum_{\l\in\s(A)}\sum_{\mu\in\s(B)}f(\l,\mu)P_\l Q_\mu,
$$
where $P_\l$ is the orthogonal projection onto $\Ker(A-\l I)$ and $Q_\mu$ is the orthogonal projection onto $\Ker(B-\mu I)$.

Similarly, one can define functions of noncommuting unitary operators. Let $f$ be a function on $\T^2$ such that \rf{absskho} holds. Then 
$$
f(\z,\t)=\sum_{j,k\in\Z}\hat f(j,k)\z^j\t^k,
$$ 
and so $f\in C(\T)\widehat\otimes C(\T)$ and
$$
f(U,V)=\sum_{j,k\in\Z}\hat f(j,k)U^jV^k
$$
for unitary operators $U$ and $V$.

It is also clear from above that functions $f$ in $B_{\be,1}^1(\T^2)$ satisfy
\rf{absskho}, and so one can define functions $f(U,V)$ for unitary $U$ and $V$ and for $f\in B_{\be,1}^1(\T^2)$.

Let us proceed now to the case of functions of triples of not necessarily commuting self-adjoint operators.

For not necessarily commuting bounded self-adjoint operators $A$, $B$ and $C$ the functions $f(A,B,C)$ can be defined as triple operator integrals
$$
f(A,B,C)=\iiint_{\R^3}f(x,y,z)\,dE_A(x)\,dE_B(y)\,dE_C(z).
$$
It was shown in \cite{Pe3}
that triple operator integrals were defined in the case when the integrand $f$ belongs to the {\it integral projective tensor product} 
$L^\be_{E_1}\widehat\otimes_{\rm i}L^\be_{E_2}\widehat\otimes_{\rm i}L^\be_{E_3}$.
Later triple operator integrals were defined in \cite{JTT} for functions in the Haagerup tensor product 
$L^\be_{E_1}\!\otimes_{\rm h}\!L^\be_{E_2}\!\otimes_{\rm h}\!L^\be_{E_3}$, see also 
\cite{AP2} for properties of such triple operator integrals.

As in the case of functions of two noncommuting operators we can define functions
$f(A,B,C)$ for functions $f$ whose Fourier transform is integrable:
$$
\int_{\R^3}|(\F f)(s,t,u)|\,ds\,dt\,du<\be
$$
In this case
$$
f(A,B,C)=\iiint_{\R^3}(\F f)(s,t,u)e^{{\rm i}sA}e^{{\rm i}tB}e^{{\rm i}uC}
\,ds\,dt\,du.
$$
However, using results of \cite{Ki} (see also \cite{N} and \cite{dLKK}), one can show that unlike the case of two operators, a functions in $B_{\be,1}^1(\R^3)$ does not have to coincide locally with a function with summable Fourier transform.

On the other hand by analogy with \rf{absolyutno!} one can show that functions in the Besov space $B_{\be,1}^2(\R^3)$ coincide locally with functions with absolutely convergent Fourier transform. 

As in the case of functions of two operators we can represent an arbitrary $2\pi$-periodic function of class $B_{\be,1}^2$  as an element of the projective tensor product 
$C(\T)\widehat\otimes C(\T)\widehat\otimes C(\T)$ in the following way:
$$
f(s,t,u)=\sum_{n\ge1}\sum_{j=-2^n}^{2^n}e^{{\rm i}js}
\sum_{k=-2^n}^{2^n}e^{{\rm i}kt}
\left(\sum_{l=-2^n}^{2^n}\widehat f_{n-1}(j,k,l)e^{{\rm i}lu}
\right).
$$
This implies that
$$
\|f\|_{C(\T)\widehat\otimes C(\T)\widehat\otimes C(\T)}\le
\sum_{n\ge1}\big(2^{n+1}+1\big)^2\|f_{n-1}\|_{L^\be(\T)}
\le\const\|f\|_{B_{\be,1}^2}.
$$
Under the assumptions that $\|A\|<\pi$, $\|B\|<\pi$ and $\|C\|<\pi$,
we can define now $f(A,B,C)$ by
$$
f(A,B,C)=\sum_{n\ge1}\sum_{j=-2^n}^{2^n}e^{{\rm i}jA}
\sum_{k=-2^n}^{2^n}e^{{\rm i}kB}
\left(\sum_{l=-2^n}^{2^n}\widehat f_{n-1}(j,k,l)e^{{\rm i}lC}
\right).
$$
In the case when $A$, $B$ and $C$ have finite spectra, we can define $f(A,B,C)$ 
for arbitrary functions on $\s(A)\times\s(B)\times\s(C)$ by
$$
f(A,B,C)=\sum_{\l\in\s(A)}\sum_{\mu\in\s(B)}\sum_{\nu\in\s(C)}
f(\l,\mu,\nu)P_\l Q_\mu R_\nu,
$$
where $P_\l$ is the orthogonal projection onto $\Ker(A-\l I)$, $Q_\mu$ is the orthogonal projection onto $\Ker(B-\mu I)$ and $R_\nu$ is the orthogonal projection onto $\Ker(C-\nu I)$.

\

\section{\bf The case of self-adjoint operators}
\label{samosop2}

\

The main purpose of this section is for each $p>2$, to construct a function $f$ of class
$B^1_{\be,1}(\R^2)$ and pairs of bounded noncommuting self-adjoint operators
$(A_1,B_1)$, $(A_2,B_2)$ such that $A_2-A_1\in\bS_p$, $B_2-B_1\in\bS_p$ but
$f(A_2,B_2)-f(A_1,B_1)\not\in\bS_p$. Then we impose additional assumptions on the support of the Fourier transform of $f$.

Denote by $\mathscr E_\s^\be(\R^d)$, $\s>0$, the space of bounded (continuous) functions $f$ on $\R^d$ whose 
Fourier transform $\F f$ is supported in 
$[-\s,\s]^d$. It is well known that 
 $\mathscr E_\s^\be(\R^2)$ equipped with the $L^\be$-norm is a Banach space.
It is easy to see that $f\in\mathscr E_\1^\be(\R^d)$ if and only if  $f(\s x)\in\mathscr E_\s^\be(\R^d)$.

It can be seen from the proof of Theorem 8.1 in \cite{ANP}  that the following result holds:

\begin{lem}
\label{2N2N+1}
For each positive integer $N$, there exist a function $f$ of class $\E_{2\pi}^\be(\R^2)$ with $\|f\|_{L^\be(\R^2)}=1$
and positive self-adjoint operators $A_1$, $A_2$ and $B$ on $\C^N$ such that $\|A_1\|=2N$,
$\|A_2\|=2N+1$, $\|B\|=N$, $f(A_2,B)=\0$  and 
$$
\|f(A_1,B)-f(A_2,B)\|_{\bS_p}=\|f(A_1,B)\|_{\bS_p}
>\const N^{\frac12-\frac1p}\|A_1-A_2\|_{\bS_p}.
$$
\end{lem}

\begin{cor} 
\label{ENAB2}
For each positive integer $N$, there exist a function $f$ of class $\E_{N}^\be(\R^2)$ with $\|f\|_{L^\be(\R^2)}=1$
and (positive) self-adjoint contractions $A_1$, $A_2$ and $B$ on a finite-dimensional
Hilbert space such that $f(A_2,B)=\0$  and 
$$
\|f(A_1,B)-f(A_2,B)\|_{\bS_p}=\|f(A_1,B)\|_{\bS_p}
>\const N^{\frac32-\frac1p}\|A_1-A_2\|_{\bS_p}.
$$
\end{cor}

\Pf Let us fix a positive integer $n$. Suppose that $f_0$, $A_1^{(0)}$, $A_2^{(0)}$ and $B^{(0)}$ satisfy the requirements of Lemma \ref{2N2N+1} with $N=n$. Put 
$f(x,y)=f_0((2n+1)x,(2n+1)y)$, $A_1=(2n+1)^{-1}A_1^{(0)}$, $A_2=(2n+1)^{-1}A_2^{(0)}$ and $B=(2n+1)^{-1}B^{(0)}$.
Then $\|f\|_{L^\be(\R^2)}=1$, $f\in\mathscr E_{2\pi(2n+1)}^\be(\R^2)\subset\mathscr E_{20n}^\be(\R^2)$,
$A_1$, $A_2$ and $B$ are self-adjoint contractions on a finite-dimensional
Hilbert space, $f(A_2,B)=\0$  and 
$$
\|f(A_1,B)-f(A_2,B)\|_{\bS_p}>\const n^{\frac12-\frac1p}(2n+1)\|A_1-A_2\|_{\bS_p}>\const(20n)^{\frac32-\frac1p}\|A_1-A_2\|_{\bS_p}.
$$
Thus, we have proved Corollary \ref{ENAB2} for $N=20n$. The general case can easily be reduced to this special case.
Let $N$ be an arbitary positive integer. I suffices to consider the case where $N>20$. Then $N=20n+k$
with $n\ge1$ and $k$, $0\le k<20$. It remains to observe that $20n>\const N$. $\bl$

\begin{cor} 
\label{ENAB3}
For each positive integer $N$, there exist a function $f$ of class $B_{\be,1}^1(\R^2)$ with $\|f\|_{B_{\be,1}^1(\R^2)}=1$
and (positive) self-adjoint contractions $A_1$, $A_2$ and $B$ on a finite-dimensional
Hilbert space such that $f(A_2,B)=\0$  and 
$$
\|f(A_1,B)-f(A_2,B)\|_{\bS_p}=\|f(A_1,B)\|_{\bS_p}
>\const N^{\frac12-\frac1p}\|A_1-A_2\|_{\bS_p}.
$$
\end{cor}

\Pf It suffices to observe that $\|f\|_{B_{\be,1}^1(\R^2)}\le\const N\|f\|_{L^\be(\R^2)}$ for all 
$f\in\mathscr E_N^\be(\R^2)$. $\bl$

Put $B_{\be,1}^1([-1,1]^2)\df\{f\big|[-1,1]^2:~f\in B_{\be,1}^1(\R^2)\}$ and 
$$
\|g\|_{B_{\be,1}^1([-1,1]^2)}\df\inf\{\|f\|_{B_{\be,1}^1(\R^2)}:~f\in B_{\be,1}^1(\R^2),~f\big|[-1,1]^2=g\}
$$
for $g\in B_{\be,1}^1([-1,1]^2)$.

Let $p\in[1,\be]$. Denote by
$\BL_p([-1,1]\times[-1,1])$ the set of functions $f$ in $B_{\be,1}^1([-1,1]^2)$ such that there exists a positive number $c$, for which
\bay
\label{AjBj}
\|f(A_1,B_1)-f(A_2,B_2)\|_{\bS_p}\le c(\|A_1-A_2\|_{\bS_p}+\|B_1-B_2\|_{\bS_p})
\ey
whenever $A_1$, $A_2$, $B_1$ and $B_2$ are
arbitrary self-adjoint contractions on a finite-dimensional
Hilbert space. We denote by $\|f\|_{\OL_{\bS_p}([-1,1]\times[-1,1])}$ the smallest constant $c$,
for which \rf{AjBj} holds. 


It is easy to see that $f\in\BL_p([-1,1]\times[-1,1])$ if and only if there exists a constant $c$ such that
$$
\|f(A_1,B)-f(A_2,B)\|_{\bS_p}\le c\|A_1-A_2\|_{\bS_p}\quad\text{and}\quad
\|f(A,B_1)-f(A,B_2)\|_{\bS_p}\le c\|B_1-B_2\|_{\bS_p}
$$
for all self-adjoint contractions $A$, $A_1$, $A_2$, $B$, $B_1$ and $B_2$ on a finite-dimensional
Hilbert space. Put 
$$
\|f\|_{\BL_p}\df\max\big\{\|f\|_{B^1_{\infty,1}([-1,1]^2)},\|f\|_{\OL_{\bS_p}([-1,1]\times[-1,1])}\big\}.
$$

Clearly, $\{f\in\BL_p([-1,1]\times[-1,1]): f(0,0)=0\}$ equipped with norm $\|\cdot\|_{\BL_p}$ is a Banach space.

\begin{lem} 
\label{Bbe11OLpSA}
Let $p\in(2,+\infty]$. Then there exists a function $f$ in $B_{\be,1}^1(\R^2)$ such that $f\big|[-1,1]^2\not\in\BL_p([-1,1]\times[-1,1])$.
\end{lem}

\Pf Suppose that $f\big|[-1,1]^2\in\BL_p([-1,1]\times[-1,1])$ for all $f\in B_{\be,1}^1(\R^2)$. 
Then by the closed graph theorem there exists a constant $c$
such that $\|f\|_{\BL_p([-1,1]\times[-1,1])}\le c\|f\|_{B^1_{\infty,1}([-1,1]^2)}$, and we get a contradiction with
Corollary \ref{ENAB3}. $\bl$

\begin{cor} 
\label{Bbe11OLpcorSA}
Let $p\in(2,+\infty]$. Then there exists a function $f$ in $B_{\be,1}^1(\R^2)$ such that for every $M>0$,
there exist self-adjoint contractions $A_1$, $A_2$ and $B$ on a finite-dimensional
Hilbert space such that
$$
\|f(A_1,B)-f(A_2,B)\|_{\bS_p}>M\|A_1-A_2\|_{\bS_p}.
$$
\end{cor}

\begin{thm} 
\label{zam8100sa}
Let $p\in(2,\be]$. Then there exists a function $f$ in $B_{\be,1}^1(\R^2)$ such that 
for any positive numbers $M$ and $\d$, there exist
self-adjoint contractions 
$A_1$, $A_2$ and $B$  on a finite-dimensional
Hilbert space such that
\bay
\label{2nervasa}
\|A_1-A_2\|_{\bS_p}<\d\quad\mbox{and}\quad
\|f(A_1,B)-f(A_2,B)\|_{\bS_p}>M\|f\|_{B_{\be,1}^1(\R^2)}\|A_1-A_2\|_{\bS_p}.
\ey
\end{thm}

\Pf 
It follows from Corollary \ref{Bbe11OLpcorSA} that  there are self-adjoint contractions $A_1$, $A_2$ and $B$  
on a finite-dimensional Hilbert space that satisfy the right inequality in \rf{2nervasa}.

We can select a positive integer $N$ such that
$\frac1N\|A_1-A_2\|_{\bS_p}<\d$. Put now $A^{(j)}=A_1+\frac jN(A_2-A_1)$, $0\le j<N$.
It remains to observe that 
$$
\big\|f(A^{(j)},B)-f(A^{(j+1)},B)\big\|_{\bS_p}>
M\|f\|_{L^\be(\R^2)}\big\|A^{(j)}-A^{(j+1)}\big\|_{\bS_p}
$$
at least for one $j$, $0\le j<N$.
$\bl$

\begin{thm}
\label{ne vkhodit00sa} 
For each $p$ in $(2,\be)$, there exist a function $f$ in 
$B^1_{\infty,1}(\R^2)$ 
and self-adjoint contractions $A_1$, $A_2$ and $B$ on Hilbert space such that 
$A_1-A_2\in\bS_p$ but $f(A_1,B)-f(A_2,B)\not\in\bS_p$.
\end{thm}

\Pf
By Theorem \ref{zam8100sa}, for each positive integer $k$, there exist self-adjoint contractions
$A_1^{(k)}$, $A_2^{(k)}$ and $B^{(k)}$ such that
$$
\big\|f\big(A_1^{(k)},B^{(k)}\big)-f\big(A_2^{(k)},B^{(k)}\big)\big\|_{\bS_p}
>1\quad{\text and}\quad\big\|A_1^{(k)}-A_2^{(k)}\big\|_{\bS_p}<2^{-k}
$$
for every $k$. It remains to put $A_1=\bigoplus\limits_{k=1}^\be A_1^{(k)}$,
$A_2=\bigoplus\limits_{k=1}^\be A_2^{(k)}$ and $B=\bigoplus\limits_{k=1}^\be B^{(k)}$. $\bl$

Let $\L$ be a closed subset of $\R^2$. Denote by $\big(B^1_{\be,1}\big)_\L\big(\R^2\big)$ the set of all functions
$f$ in $B^1_{\be,1}\big(\R^2\big)$ such that $\supp\mathscr F f\subset\L$.
The seminorm $\|\cdot\|_{B^1_{\be,1}(\R^2)}$ is a norm on $\big(B^1_{\be,1}\big)_\L\big(\R^2\big)$ if and only if $0\notin\L$.


Let $\L$ be a closed subset of $\R^2$ and $r>0$ . Put $\t_\L(r)\df\inf R$, where the infimum is taken over all $R>r$
such that the set $\L\cap R\dd$ contains a disk of radius $r$. In particular, $\t_\L(r)=\be$ if 
$\L$ contains no disk of radius $r$.

\

\begin{thm}
\label{uniszhaLSA}
Let $\L\subset\R^2$ and $2<p<\be$. Suppose that $\liminf\limits_{r\to\be}r^{\frac1p-\frac32}\t_\L(r)=0$. 
Then there exists a function $f$ in $(B^1_{\be,1})_\L(\R^2)$ and self-adjoint operators 
$A_1$, $A_2$ and $B$ such that $A_1-A_2\in\bS_p$ but 
$f(A_1,B)-f(A_2,B)\not\in\bS_p$.
\end{thm}

We need the following lemma.

\begin{lem} 
\label{29}
Let $R>r\ge1$ and let $D$ be a disk of radius $r$, $D\subset R\dd$. 
Then there exist a function $f\in\mathscr E_{R}^\be(\R^2)$ with $\supp\mathscr F f\subset D$
and positive self-adjoint contractions $A_1$, $A_2$ and $B$ on a finite-dimensional
Hilbert space such that 
$$
\|f(A_1,B)-f(A_2,B)\|_{\bS_p}>\const r^{\frac32-\frac1p}R^{-1}\|f\|_{B^1_{\be,1}(\R^2)}\|A_1-A_2\|_{\bS_p}.
$$
\end{lem}
\Pf It suffices to consider the case when $r$ is a positive integer. Applying Corollary \ref{ENAB2}
to $N=r$, we get a function $f_0\in\mathscr E_{r}^\be(\R^2)$ with $\|f\|_{L^\be(\R^2)}=1$ and 
positive self-adjoint contractions $A_1$, $A_2$ and $B$ on a finite-dimensional
Hilbert space such that $f_0(A_2,B)=\0$ and
$$
\|f_0(A_1,B)-f_0(A_2,B)\|_{\bS_p}>\const r^{\frac32-\frac1p}\|A_1-A_2\|_{\bS_p}.
$$
We can take $a\in\R^2$ such that $|a|<R$ in such a way that $a+r\dd\subset D$. 
Put $f\df e^{{\rm i}(x,a)}f_0$. Then
$$
\|f(A_1,B)-f(A_2,B)\|_{\bS_p}=\|f_0(A_1,B)-f_0(A_2,B)\|_{\bS_p}>\const r^{\frac32-\frac1p}\|A_1-A_2\|_{\bS_p}.
$$
It remains to observe that $\|f\|_{B^1_{\be,1}(\R^2)}\le\const R\|f\|_{L^\infty(\R^2)}=\const R$. $\bl$

\begin{lem} 
\label{Bbe11OLpLa}
Let $\L\subset\R^2$ and let $2<p<\be$. Suppose that $\liminf\limits_{r\to\be}r^{\frac1p-\frac32}\varkappa_\L(r)=0$. 
Then there exists a function $f$ in $(B_{\be,1}^1)_\L(\R^2)$ such that $f\big|[-1,1]^2\not\in\BL_p([-1,1]\times[-1,1])$.
\end{lem}

\Pf To prove the result, we can repeat the the argument in the proof of Lemma \ref{Bbe11OLpSA} and use Lemma \ref{29}
instead of Corollary \ref{ENAB3}. $\bl$

\medskip

This lemma implies Theorem  \ref{uniszhaLSA} in the same a way as Lemma \ref{Bbe11OLpSA}
implies Theorem \ref{ne vkhodit00sa}.

Theorem \ref{uniszhaLSA} readily implies the following theorem.

\begin{thm}
\label{uniszhaLcorSA}
Let $\L\subset\R^2$ and $1\le\a<\frac32$. Suppose that $\t_\L(r)\le\const r^\a$ for all sufficiently large $r$. Then for each 
$p>(\frac32-\a)^{-1}$, there exists a function $f$ in $(B^1_{\be,1})_\L(\R^2)$ and self-adjoint operatrs  
$A_1$, $A_2$ and $B$ such that 
$A_1-A_2\in\bS_p$ but $f(A_1,B)-f(A_2,B)\not\in\bS_p$.
\end{thm}

\begin{cor}
\label{uniszhaLcrSA}
Let $\L$ be a nondegenerate angle in $\R^2$.
 Then for each 
$p>2$, there exists a function $f$ in $(B^1_{\be,1})_\L(\R^2)$ and self adjoint operators  $A_1$, $A_2$ and $B$ such that 
$A_1-A_2\in\bS_p$ but $f(A_1,B)-f(A_2,B)\not\in\bS_p$.
\end{cor}

The following theorem is a version of Theorem \ref{ne vkhodit00sa} in the case $p=\infty$.

\begin{thm}
\label{ne vkhodit00sabe} 
There exist a function $f$ in 
$B^1_{\infty,1}(\R^2)$ 
and self-adjoint contractions $A_1$, $A_2$ and $B$ on Hilbert space such that 
$A_1-A_2$ is compact but $f(A_1,B)-f(A_2,B)$ is not compact.
\end{thm}

\Pf
By Theorem \ref{zam8100sa}, for each positive integer $k$, there exist self-adjoint contractions
$A_1^{(k)}$, $A_2^{(k)}$ and $B^{(k)}$ on a finite-dimensional
Hilbert space $\mathcal H_k$ such that
such that
$$
\big\|f\big(A_1^{(k)},B^{(k)}\big)-f\big(A_2^{(k)},B^{(k)}\big)\big\|>
1\quad{\text and}\quad\big\|A_1^{(k)}-A_2^{(k)}\big\|<2^{-k}
$$
for every $k$. It remains to put $A_1=\bigoplus\limits_{k=1}^\be A_1^{(k)}$,
$A_2=\bigoplus\limits_{k=1}^\be A_2^{(k)}$, $B=\bigoplus\limits_{k=1}^\be B^{(k)}$
and observe that $A_1-A_2$ is compact but $f(A_1,B)-f(A_2,B)$ is not compact. $\bl$

\medskip

In the same way we can get  a version of Theorem \ref{uniszhaLSA} in the case $p=\infty$.
The same is true about Theorem \ref{uniszhaLcorSA} and Corollary \ref{uniszhaLcrSA}, which
are special cases of Theorem \ref{uniszhaLSA}.

\

\section{\bf The case of contractions}
\label{szhali!}

\

The purpose of this section is to obtain analogs of the results of \S\:\ref{samosop2} for functions of pairs of noncommuting contractions.
We obtain analogs of Theorems \ref{ne vkhodit00sa}, \ref{uniszhaLSA}, \ref{uniszhaLcorSA} and \ref{ne vkhodit00sabe}.

Actually, we are going to work with unitary operators and obtain analogs of Theorems \ref{ne vkhodit00sa}, 
and \ref{ne vkhodit00sabe} for unitary operators and functions that belong to the Besov class 
$\Ba$ of functions analytic in the bidisk. Since unitary operators are contractions, we obtain thereby analogs of 
Theorems \ref{ne vkhodit00sa} and \ref{ne vkhodit00sabe}  for functions of contractions.

%
 
Let $p\in[1,\be]$. Denote by
$\BL_p(\T\times\T)$ the set of the functions $f$ of class $\Ba$ such that
\bay
\label{UjVj}
\|f(U_1,V_1)-f(U_2,V_2)\|_{\bS_p}\le c(\|U_1-U_2\|_{\bS_p}+\|V_1-V_2\|_{\bS_p})
\ey
for arbitrary unitary operators $U_1$, $U_2$, $V_1$ and $V_2$ on a finite-dimensional
Hilbert space. We denote by $\|f\|_{\OL_{\bS_p}(\T\times\T)}$ the smallest constant $c$,
for which \rf{UjVj} holds.


It is easy to see that $f\in\BL_p(\T\times\T)$ if and only if there exists a constant $c$ such that
$$
\|f(U_1,V)-f(U_2,V)\|_{\bS_p}\le c\|U_1-U_2\|_{\bS_p}\quad\text{and}\quad
\|f(U,V_1)-f(U,V_2)\|_{\bS_p}\le c\|V_1-V_2\|_{\bS_p}
$$
for all unitary operators $U$, $U_1$, $U_2$, $V$, $V_1$ and $V_2$ on a finite-dimensional
Hilbert space.
Put 
$$
\|f\|_{\BL_p}\df\max\big\{\|f\|_{B^1_{\infty,1}(\T^2)},\|f\|_{\OL_{\bS_p}(\T\times\T)}\big\}.
$$
Clearly, the space $\{f\in\BL_p(\T\times\T): f(1,1)=0\}$ equipped with the norm $\|\cdot\|_{\BL_p}$ is a Banach space.

\begin{lem} 
\label{Bbe11OLp}
Let $p\in(2,+\infty]$. Then there exists a function $f$ of class $\Ba$ such that $f\not\in\BL_p(\T\times\T)$.
\end{lem}

\Pf Suppose that $\Ba\subset\BL_p(\T\times\T)$. It is easy to deduce from the closed graph theorem there exists a constant $c$
such that $\|f\|_{\BL_p(\T\times\T)}\le c\|f\|_{B^1_{\infty,1}(\T^2)}$. This contradicts 
Theorem 7.3 of \cite{AP}. $\bl$

\begin{cor} 
\label{Bbe11OLpcor}
Let $p\in(2,+\infty]$. Then there exists a function $f\in\Ba$ such that for every $M>0$, 
there exist unitary operators $U_1$, $U_2$ and $V$ on a finite-dimensional
Hilbert space such that
$$
\|f(U_1,V)-f(U_2,V)\|_{\bS_p}>M\|U_1-U_2\|_{\bS_p}.
$$
\end{cor}

\medskip

\begin{lem} 
\label{UASp}
Let $U$ be a unitary operator such that $U-I\in\bS_p$, where $p\in[1,\be]$.
Then there exists a self-adjoint operator $A$ such that $e^{{\rm i}A}=U$
and
$$
\|U-I\|_{\bS_p}\le\|A\|_{\bS_p}\le\frac\pi2\|U-I\|_{\bS_p}.
$$
\end{lem}

\Pf We can take a self-adjoint operator $A$ such that
$e^{{\rm i}A}=U$ and $\|A\|\le\pi$. Clearly, $\|e^{{\rm i}A}-I\|=|e^{{\rm i}\|A\|}-1|$. That implies the required
result for $p=\be$.
Now let $p<\be$.  Clearly, $A$ is a compact self-adjoint operator. Let $\{s_k\}_{k=1}^\be$ be
a sequence of all eigenvalues taking in account the multiplicity. Then
$$
\|U-I\|_{\bS_p}^p=\sum_{k=1}^\be|e^{{\rm i}s_k}-1|^p\le\sum_{k=1}^\be|s_k|^p=\|A\|_{\bS_p}^p
\le\left(\frac\pi2\right)^p\sum_{k=1}^\be|e^{{\rm i}s_k}-1|^p=\left(\frac\pi2\right)^p\|U-I\|_{\bS_p}^p.\, \bl
$$

\begin{cor} 
\label{UASpcor}
Let $1\le p\le\be$ and
let $U_1$ and $U_2$ be unitary operators such that $U_1-U_2\in\bS_p$.
Then for each positive integer $N$ there exists a sequence 
$\big\{U^{[k]}\big\}_{k=0}^N$ of unitary operators
such that $U^{[0]}=U_1$, $U^{[N]}=U_2$ and 
$\big\|U^{[k]}-U^{[k-1]}\big\|_{\bS_p}\le\dfrac\pi{2N}\|U_1-U_2\|_{\bS_p}$ for 
$k=1,2,\cdots, N$.
\end{cor}

\Pf By Lemma \ref{UASp}, there exists a self-adjoint operator $A$ such that $e^{{\rm i}A}=U_1^{-1}U_2$
and $\|A\|_{\bS_p}\le\dfrac\pi2\|U_2U_1^{-1}-I\|_{\bS_p}=\dfrac\pi2\|U_1-U_2\|_{\bS_p}$.
It remains to put $U^{[k]}=U_1e^{{\rm i}\frac kNA}$. $\bl$

\medskip


\begin{thm}
\label{UdM} 
Suppose that $2<p\le\be$. Then there exists a function $f$ in $\Ba$ such that for any positive numbers $M$ and $\d$, there exist a finite-dimensional Hilbert space and  unitary operators  $U_1$, $U_2$ and $V$ on it such that 
$$
\|U_1-U_2\|_{\bS_p}<\d\qquad\mbox{but}\qquad
\|f(U_1,V)-f(U_2,V)\|_{\bS_p}>M\|U_1-U_2\|_{\bS_p}.
$$
\end{thm} 

\Pf By Corollary \ref{Bbe11OLpcor} there exists a function  $f\in\Ba$ such that for every $M>0$,
there exist unitary operators $U_1$, $U_2$ and $V$ such that
$$
\|f(U_1,V)-f(U_2,V)\|_{\bS_p}>\frac\pi2M\|U_1-U_2\|_{\bS_p}.
$$
We can take an integer $N$ such that $\dfrac\pi{2N}<\d$. Applying Corollary \ref{UASpcor},
we get a sequence $\{U^{[k]}\}_{k=0}^N$ of unitary operators such that 
$U^{[0]}=U_1$, $U^{[N]}=U_2$ and 
$$
\big\|U^{[k]}-U^{[k-1]}\big\|_{\bS_p}\le\dfrac\pi{2N}\|U_1-U_2\|_{\bS_p}<\d\|U_1-U_2\|_{\bS_p},\quad k=1,2,\cdots,N.
$$ 
It remains to prove that 
$\big\|f\big(U^{[k]},V\big)-f\big(U^{[k-1]},V\big)\big\|_{\bS_p}>
M\big\|U^{[k]}-U^{[k-1]}\big\|_{\bS_p}$
for some integer $k$ such that $1\le k\le N$.

Assume the contrary. Then
\begin{align*}
\|f(U_1,V)-f(U_2,V)\|_{\bS_p}&\le
\sum_{k=1}^N\big\|f\big(U^{[k]},V)-f\big(U^{[k-1]},V\big)\big\|_{\bS_p}\\
&\le M\sum_{k=1}^N\big\|U^{[k]}-U^{[k-1]}\big\|_{\bS_p}
\le\frac\pi2M\|U_1-U_2\|_{\bS_p},
\end{align*}
and we get a contradiction. $\bl$

\medskip


\begin{thm}
\label{uniszha}
Suppose that $2<p<\be$. Then there exists a function $f$ in $\Ba$ and unitary operators  $U_1$, $U_2$ and $V$ such that $U_1-U_2\in\bS_p$ but $f(U_1,V)-f(U_2,V)\not\in\bS_p$.
\end{thm}

\Pf
By Theorem \ref{UdM}, for each positive integer $k$, there exist unitary operators
$U_{1,k}$, $U_{2,k}$ and $V_k$ such that
$$
\big\|f\big(U_{1,k},V_k\big)-f\big(U_{2,k},V_k\big)\big\|_{\bS_p}
>1\quad\text {and}\quad\big\|U_{1,k}-U_{2,k}\big\|_{\bS_p}<2^{-k}
$$
for every $k$. It remains to put $U_1=\bigoplus\limits_{k=1}^\be U_{1,k}$,
$U_2=\bigoplus\limits_{k=1}^\be U_{2,k}$ and $V=\bigoplus\limits_{k=1}^\be V_k$. 
$\bl$

Let $\L$ be a subset of $\Z^2$. Put $(B^1_{\be,1})_\L(\T^2)\df\{f\in B^1_{\be,1}(\T^2): \supp\widehat f\subset\L\}$.
For a positive integer $m$ we, denote by $\varkappa_\L(m)$
the least positive integer $N$ such that there exist $n_1, n_2\in\Z$ satisfying the following conditions
$$
-N\le n_1,n_2\le N-m\quad\mbox{and}
\quad([n_1,n_1+m]\times[n_2,n_2+m])\cap\Z^2\subset\Lambda.
$$
Put $\varkappa_\L(m)\df\be$ if there is no such $N$.

The same technique as above allows us to obtain the following result:

\begin{thm}
\label{uniszhaL}
Let $\L\subset\Z^2$ and $2<p<\be$. Suppose that $\liminf\limits_{m\to\be}m^{\frac1p-\frac32}\varkappa_\L(m)=0$. 
Then there exists a function $f$ in $(B^1_{\be,1})_\L(\T^2)$ and unitary operators  
$U_1$, $U_2$ and $V$ such that $U_1-U_2\in\bS_p$ but 
$f(U_1,V)-f(U_2,V)\not\in\bS_p$.
\end{thm}

To prove this theorem, we need the following result which is essentially Lemma 7.2 in \cite{AP}; see also its proof in \cite{AP}.

\begin{lem} 
\label{fU1U2V}
For each $m\in\Bbb N$, there exists an analytic polynomial $f$ in two variables of degree at most $4m-2$ in each variable, and unitary operators $U_1$, $U_2$
and $V$ in $\C^m$ such that 
\bay
\label{fU2V0}
f(U_2,V)=\0\quad\text{and}\quad\|f(U_1,V)
\|_{\bS_p}>\pi^{-1}m^{\frac32-\frac1p} \|f\|_{L^\be(\T^2)}\|U_1-U_2\|_{\bS_p}
\ey
for every $p>0$.
\end{lem}

\medskip

{\bf Remark.}  Put $g(z_1,z_2)=z_1^{n_1}z_2^{n_2}f(z_1,z_2)$, where $n_1,n_2\in\Z$. If $f$ satisfies \rf{fU2V0},
then $g$ also satisfies \rf{fU2V0}.

\medskip

This remark implies the following result:

\begin{lem} 
\label{fU1U2VSQ}
Let $\L\subset\Z^2$ and let $2<p<\be$. Suppose that 
$$
\Lambda\supset([n_1,n_1+m]\times[n_2,n_2+m])\cap\Z^2
$$
for some integers $n_1$, $n_2$ such that 
$-N\le n_1,n_2\le N-m$ and 
$$
([n_1,n_1+m]\times[n_2,n_2+m])\cap\Z^2\subset\Lambda
$$
for some positive integers $m$ и $N$.
Then there exist a trigonometric polynomial $f$ in $(B^1_{\be,1})_\L(\T^2)$ and unitary operators $U_1$, $U_2$ and $V$ on a finite-dimensional
Hilbert space such that
$$
\|f(U_1,V)-f(U_2,V)\|_{\bS_p}>\const m^{\frac32-\frac1p} N^{-1}\|f\|_{B^1_{\be,1}}\|U_1-U_2\|_{\bS_p}.
$$
\end{lem}

{\bf Sketch of the proof of Lemma \ref{fU1U2VSQ}.} It suffices consider the case where $m\in16\Bbb N$. Let $m=16m_0$.
Note that there exists a ``square''  
$$
Q'=([n_1',n_1'+4m_0]\times[n_2',n_2'+4m_0])\cap\Z^2
$$
that is contained in the ``square'' 
$$
([n_1,n_1+m]\times[n_2,n_2+m])\cap\Z^2
$$ 
and such that 
$\dist((0,0),Q')\le\const m$. It remains to apply the remark after Lemma \ref{fU1U2V} and to observe
that $\|g\|_{B^1_{\be,1}}\le \const N\|g\|_{L^\be(\T^2)}$ for any trigonometrical polynomial $g$ with
$\supp\widehat g\subset Q'$. $\bl$

\begin{cor}
\label{000}
Under the hypotheses of Theorem {\em\ref{uniszhaL}}, for every $M>0$, there exist a function $f\in(B^1_{\be,1})_\L(\T^2)$
and unitary operators 
$U_1$, $U_2$ and $V$ on a finite-dimensional
Hilbert space such that
$$
\|f(U_1,V)-f(U_2,V)\|_{\bS_p}>M\|f\|_{B^1_{\be,1}}\|U_1-U_2\|_{\bS_p}.
$$
\end{cor}

Theorem  \ref{uniszhaL} can be deduced now from Corollary \ref{000} in the same way as 
Theorem \ref{uniszha} was deduced from
Theorem 7.3 in \cite{AP}.

Theorem \ref{uniszhaL} readily implies the following result:

\begin{thm}
\label{uniszhaLcor}
Let $\L\subset\Z^2$ and $1\le\a<\frac32$. Suppose that $\varkappa_\L(m)\le\const m^\a$ for all $m\ge1$. Then for each 
$p>(\frac32-\a)^{-1}$, there exists a function $f$ in $(B^1_{\be,1})_\L(\T^2)$ and unitary operators  $U_1$, $U_2$ and $V$ such that 
$U_1-U_2\in\bS_p$ but $f(U_1,V)-f(U_2,V)\not\in\bS_p$.
\end{thm}

\begin{cor}
\label{uniszhaLcr}
Let $X$ be a nondegenerate angle in $\R^2$.
Let $\L\subset\Z^2$. Suppose that $\Lambda\supset X\cap\Z^2$. Then for each 
$p>2$, there exist a function $f$ in $(B^1_{\be,1})_\L(\T^2)$ and unitary operators  $U_1$, $U_2$ and $V$ such that 
$U_1-U_2\in\bS_p$ but $f(U_1,V)-f(U_2,V)\not\in\bS_p$.
\end{cor}

In the same way as in the case of self-adjoint operators (see \S\:\ref{samosop2}) all main results 
of this section have natural versions  in the case $p=\infty$. For example, 
the corresponding version of Theorem \ref{uniszha}  in the case $p=\be$ can be formulated
as follows.

\begin{thm}
\label{uniszhaco}
There exist a function $f$ in $\Ba$ and unitary operators  $U_1$, $U_2$ and $V$ such that $U_1-U_2$ is compact 
but $f(U_1,V)-f(U_2,V)$ is not compact.
\end{thm}

\

\section{\bf The case of triples of noncommuting self-adjoint operators}
\label{triple}

\

A slight modification of the proof of Theorem 4.1 of \cite{Pe4} shows that if 
$1\le p\le\be$, 
for $p\in[1,\be]$, there exist
a sequence $\{f_N\}$ of functions in $\E_{2\pi}^\be(\R^3)$, sequences 
$\big\{A^{(N)}\big\}$, $\big\{B^{(N)}\big\}$, $\big\{C_1^{(N)}\big\}$,
and $\big\{C_2^{(N)}\big\}$ of self-adjoint operators 
of rank at most $N$ such that
$$
\|f_N\|_{L^\be}\le\const,
$$
and
$$
\big\|\big(f(A^{(N)},B^{(N)},C_1^{(N)}\big)-
\big(f(A^{(N)},B^{(N)},C_2^{(N)}\big)\big\|_{\bS_p}
\ge\const N^{1/2}\|C_1-C_2\|_{\bS_p}.
$$


Recall that we have defined in \S\:\ref{Fnko} functions $f(A,B,C)$ of arbitrary bounded self-adjoint operators $A$, $B$ and $C$ for functions of class 
$B_{\be,1}^2(\R^3)$.

The method used in \S\:\ref{samosop2} allows us to prove the following result.

\begin{thm}
\label{troikaSp}
Let $1<p<\be$. Then there exist a function $f$ in $B_{\be,1}^2(\R^3)$ and triples self-adjoint operators $(A_1,B_1,C_1)$ and $(A_2,B_2,C_2)$ such that 
$A_2-A_1\in\bS_p$, $B_2-B_1\in\bS_p$, $C_2-C_1\in\bS_p$ but
$$
f(A_2,B_2,C_2)-f(A_1,B_1,C_1)\not\in\bS_p.
$$
\end{thm}

Note that as in \S\:\ref{samosop2}, one can prove a similar result if we replace 
$\bS_p$ with the class of compact operators.

\
 
\noindent
\begin{tabular}{p{7cm}p{15cm}}
A.B. Aleksandrov & V.V. Peller \\
St.Petersburg Branch & Department of Mathematics \\
Steklov Institute of Mathematics  & Michigan State University \\
Fontanka 27, 191023 St.Petersburg & East Lansing, Michigan 48824\\
Russia&USA\\
email: alex@pdmi.ras.ru&and\\
&Peoples' Friendship University\\
& of Russia (RUDN University)\\
&6 Miklukho-Maklaya St., Moscow,\\
& 117198, Russian Federation\\
& email: peller@math.msu.edu
\end{tabular}

\end{document}